\documentclass[11pt,reqno]{amsart}
\usepackage{amssymb,amscd,amsbsy,mathrsfs}
\usepackage{mathdots}
\newcommand{\lb}{\linebreak}

\newcommand{\f}{\varphi}

\newcommand{\G}{\Gamma}
\newcommand{\D}{\Delta}

\newcommand{\F}{{\mathscr F}}

\newcommand{\K}{{\mathscr K}}

\newcommand{\cp}{{\mathscr P}}

\newcommand{\T}{{\Bbb T}}
\newcommand{\pp}{{\Bbb P}}
\newcommand{\dd}{{\Bbb D}}
\newcommand{\R}{{\Bbb R}}
\newcommand{\Z}{{\Bbb Z}}

\newcommand{\0}{{\boldsymbol{0}}}

\newcommand{\bs}{\boldsymbol}

\newcommand{\bS}{{\boldsymbol S}}

\newcommand{\rf}[1]{(\ref{#1})}

\newcommand{\df}{\stackrel{\mathrm{def}}{=}}

\newcommand{\supp}{\operatorname{supp}}

\newcommand{\const}{\operatorname{const}}

\newcommand{\eeq}{\end{equation}}
\newcommand{\beq}{\begin{equation}}
\newcommand{\bay}{\begin{eqnarray}}
\newcommand{\ba}{\begin{align*}}
\newcommand{\ea}{\end{align*}}
\newcommand{\ey}{\end{eqnarray}}
\newcommand{\bey}{\begin{eqnarray*}}
\newcommand{\eey}{\end{eqnarray*}}

\newcommand{\be}{\infty}

\newcommand{\bl}{\blacksquare}

\newcommand{\Pf}{{\bf Proof. }}

\newtheorem{thm}{\hspace{\parindent}Theorem}[section]

\newtheorem{cor}[thm]{\hspace{\parindent}Corollary}

\theoremstyle{remark}

\newtheorem*{rem*}{Remark}

\newcommand{\ri}{{\rm i}}

\newcommand\fM{\frak M}

\newcommand{\Bbbone}{{\rm{1\mathchoice{\kern-0.25em}{\kern-0.25em}{\kern-0.2em}{\kern-0.2em}I}}}

\newcommand\mB{\mathcal{B}}

\begin{document}

\numberwithin{equation}{section}

\numberwithin{equation}{section}

\title{Triangular projection on $\boldsymbol{\bS_p,~0<p<1}$, \\ as $\bs{p}$ approaches 1}
\author{A.B. Aleksandrov and V.V. Peller}
\thanks{The research on \S\:3-5  is supported by 
Russian Science Foundation [grant number 23-11-00153].
The research on the rest of the paper is supported by by a grant of the Government of the Russian Federation for the state support of scientific research, carried out under the supervision of leading scientists, agreement  075-15-2021-602.}


\

\begin{abstract}
This is a continuation of our recent paper. We continue studying
properties of the triangular projection $\cp_n$ on the space of $n\times n$ matrices.  We establish sharp estimates for the $p$-norms of $\cp_n$ as an operator on the Schatten--von Neumann class $\bS_p$
for $0<p<1$. Our estimates are uniform in $n$ and $p$ as soon as $p$ is separated away from 0.
The main result of the paper shows that for $p\in(0,1)$, the $p$-norms of
$\cp_n$ on $\bS_p$
behave as $n\to\be$ and $p\to1$ as $n^{1/p-1}\min\big\{(1-p)^{-1},\log n\big\}$. 
\end{abstract} 

\maketitle


\setcounter{section}{0}
\section{\bf Introduction}
\setcounter{equation}{0}
\label{In}

\

This paper continues studying Schur multipliers of Schatten--von Neumann classes 
$\bS_p$ for $p\in(0,1)$. We refer the reader to \cite{AP1}, \cite{AP2} and \cite{AP3} for earlier results.

Recall that for infinite matrices $A=\{a_{jk}\}_{j,k\ge0}$ and 
$B=\{b_{jk}\}_{j,k\ge0}$ the {\it Schur--Hadamard product} $A\star B$ of $A$ and $B$ is, by definition,
the matrix
$$
A\star B=\{a_{jk}b_{jk}\}_{j,k\ge0}.
$$
In the same way one can define the Schur--Hadamard product of two finite matrices of the same size.

Let $p>0$. We say that a matrix $A=\{a_{jk}\}_{j,k\ge0}$ is a {\it Schur multiplier} of the Schatten--von Neumann class $\bS_p$ if
$$
B\in\bS_p\quad\Longrightarrow A\star B\in\bS_p.
$$
We use the notation  $\fM_p$ for the space of Schur multipliers of $\bS_p$. 

We refer the reader to \cite{GK1} and \cite{BS2} for properties of Schatten-von Neumann classes $\bS_p$. It is well known that for $p\ge1$, the class $\bS_p$ is a Banach space and $\|\cdot\|_{\bS_p}$ is a norm. For $p<1$, the class $\bS_p$ is a quasi-Banach space. More precisely, for $p<1$, $\bS_p$ is a
{\it $p$-Banach space} and $\|\cdot\|_{\bS_p}$ is a {\it $p$-norm}, i.e.,
$$
\|T+R\|_{\bS_p}^p\le\|T\|_{\bS_p}^p+\|R\|_{\bS_p}^p,\quad T,~R\in\bS_p,
$$
see, for example, \cite{Pe2}, App. 1, for the proof of this triangle inequality for $p<1$.

It follows from the closed graph theorem that if $A\in\fM_p$, then
$$
\|A\|_{\fM_p}\df\sup\{\|A\star B\|_{\bS_p}:~B\in\bS_p,~\|B\|_{\bS_p}\le1\}<\be.
$$
It is easy to see that $\fM_p$ is a Banach space for $p\ge1$ and it is a $p$-Banach space for $0<p<1$.

Obviously, all finite matrices are Schur multipliers of $\bS_p$ for all $p>0$ and we can define the norm $\|\cdot\|_{\fM_p}$ (the $p$-norm $\|\cdot\|_{\fM_p}$ if $p<1$) as in the case of infinite matrices.

We denote by $\mB(\bS_p)$ the space of bounded linear transformers on $\bS_p$ and 
for ${\mathcal T}\in\mB(\bS_p)$, we put
$$
\|{\mathcal T}\|_{\mB(\bS_p)}=\sup\big\{\|{\mathcal T}T\|_{\bS_p}:~T\in\bS_p,~\|T\|_{\bS_p}\le1\big\}.
$$
Clearly, $\|\cdot\|_{\mB(\bS_p)}$ is a $p$-norm for $p\in(0,1)$.

In our previous paper \cite{AP3} we 
solved a problem by B.S. Kashin\footnote{the problem was posed by B.S. Kashin at the conference dedicated to the memory of A.A. Gonchar and A.G. Vitushkin, Moscow, 2021} on the behavior of $\|\cp_n\|_{\mB(\bS_p)}$ for large values of $n$ in the case when $p\in(0,1)$, where $\cp_n$ is the upper triangular projection on the space of $n\times n$ matrices, i.e., for $A=\{a_{jk}\}_{1\le j\le n,1\le k\le n}$, the matrix 
$\cp_nA$ is defined by
$$
\cp_nA=\chi_n\star A,
$$
where the matrix $\chi_n\df\{(\chi_n)_{jk}\}_{1\le j\le n,1\le k\le n}$ is given by
$$
(\chi_n)_{jk}=\left\{\begin{array}{ll}1,&j\le k,\\0,&j>k.\end{array}\right.
$$

Let us introduce the infinite Hankel matrix\footnote{see \S\ref{Gankeli} for the definition of Hankel matrices} $\D_n\df\{(\D_n)_{jk}\}_{j,k\ge0}$ is defined by 
$$
(\D_n)_{jk}=\left\{\begin{array}{ll}1,&j+k<n,\\0,&j+k\ge n.\end{array}\right.
$$
Clearly, $\|\chi_n\|_{\fM_p}=\|\D_n\|_{\fM_p}$. Indeed, the matrix $\chi_n$ can be obtained from the matrix $\D_n$ by deleting the zero rows and the zero columns and rotating by $\pi/2$.

The main result of \cite{AP3} shows that for $0<p<1$, there exist positive numbers $c$ and $C$ such that
$$
cn^{1/p-1}\le\|\cp_n\|_{\mB(\bS_p)}=\|\chi_n\|_{\fM_p}=\|\D_n\|_{\fM_p}\le Cn^{1/p-1}.
$$

Note that in the case when $1<p<\be$, the projections $\cp_n$, $n\ge1$, are uniformly bounded in $n$ on 
$\bS_p$. This follows immediately from the well known fact that for $p\in(1,\be)$, the {\it triangular projection} $\cp$ onto the infinite upper triangular matrices is bounded, see \cite{GK2}. 
On the other hand,  in the case when $p=1$ the projection $\cp$ 
is unbounded, see \cite{GK2}.

It is also well known that there are positive numbers $k$ and $K$ such that
\bay
\label{a_chto_esli_p=1?}
k \log(1+n)\le\|\chi_n\|_{\fM_1}\le K \log(1+n),
\ey
see e.g., Remark 6 in \S\;6 of \cite{AP3}. 

In this paper we address a question posed by D.M. Stolyarov. He asked whether it is possible to 
obtain sharp estimates for $\|\cp_n\|_{\mB(\bS_p)}$ uniform in $n$ and $p$ if $p$ is separated away from 0
that would lead to estimates \rf{a_chto_esli_p=1?} as $p$ tends to 1. We give an affirmative answer to his question.

A brief introduction in the Besov spaces $B_p^{1/p}$ will be given in \S\;\ref{Besovy}.
In \S\;\ref{otsenkitrigpolinomov} we improve earlier results and obtain estimates of certain trigonometric polynomials that will be used later. 
In \S\,\ref{Gankeli} we introduce Hankel matrices and we state a description of Hankel matrices of class $\bS_p$, $0<p<\be$, in terms of Besov classes. We also state in \S\,\ref{Gankeli} an upper estimate for the $p$-norms of Hankel matrices in the space of Schur multipliers  $\fM_p$. 

Section \ref{analDir} of this paper is devoted to obtaining sharp lower and upper estimates for the $p$-norms in $H^p$ of the analytic Dirichlet kernel $D_n^+$. In \S\:\ref{povedenie} we reduce studying the behaviour of 
$\|\cp_n\|_{\mB(\bS_p)}$ to the behaviour of $\|D_n^+\|_{H^p}$. The main result of the paper is given in 
\S\:\ref{Osnova}. Namely, we show that the $p$-norms $\|\cp_n\|_{\mB(\bS_p)}$ behave like
\lb$n^{1/p-1}\min\big\{(1-p)^{-1},\log n\big\}$ as $n\to\be$ and $p\to1$.

\

\section{\bf The Besov classes $\bs{B_p^{1/p}}$}
\setcounter{equation}{0}
\label{Besovy}

\

For the reader's convenience we give here the definition of the Besov classes 
$\big(B_p^{1/p}\big)_+$ of functions in $B_p^{1/p}$ that are analytic in the unit disk $\dd$.

\medskip

{\bf Definition.} Let $v$ be an infinitely differentiable function on $(0,\be)$ such that
$v\ge\0$, $\supp v=\big[\frac12,2\big]$ and
$$
\sum_{j\ge0}v\big(2^{-j}x\big)=1\quad\mbox{for}\quad x\ge1.
$$
Let $0<p<\be$. For a function $\f$ analytic in $\dd$ we say that 
$\f\in\big(B_p^{1/p}\big)_+$ if
$$
\sum_{n\ge0}2^n\|\f*V_n\|_{L^p}^p<\be,
$$
where 
\bay
\label{Vn}
V_n(z)=\sum_{j>0}v\big(2^{-n}j\big)z^j
\quad\mbox{and}\quad
\f*V_n=\sum_{j>0}\widehat\f(j)\widehat V_n(j)z^j.
\ey
Among various equivalent $p$-norms on $\big(B_p^{1/p}\big)_+$ we can select the following one:
$$
\|\f\|_{B_p^{1/p}}\df\left(\sum_{n\ge0}2^n\|\f*V_n\|_{L^p}^p\right)^{1/p}.
$$

\

\section{\bf  $\bs{L^p}$ estimates for certain trigonometric polynomials}
\setcounter{equation}{0}
\label{otsenkitrigpolinomov}

\

We start this section with an elementary observation. Since clearly, $|f(0)|\le\|f\|_{H^p}$, $p\in(0,\be]$ for $f\in H^p$, it follows that
$$
\mbox{if}\quad f(z)=\sum\limits_{j=m}^n\widehat f(j)z^j,\quad\mbox{then}\quad
|\widehat f(m)|\le\|f\|_{L^p}\quad\mbox{and}\quad |\widehat f(n)|\le\|f\|_{L^p}\quad p>0.
$$

\medskip


We need sharp estimates of the $p$-norms in $L^p$ of certain trigonometric polynomials whose coefficients are expressed in terms of the values of an infinitely differentiable function with compact support at the points of the grid $\{\frac km\}_{k\in\Z}$.

Let $q$ be an infinitely differentiable function on $\R$ with compact support, we define the family $Q_m$, $m\ge1$, of trigonometric polynomials by
\bay
\label{Ups}
Q_m(z)=\sum_{k\in\Z}q\left(\frac km\right)z^k,\quad m\ge1.
\ey
It was established in \cite{A}, Ch 4 that for $p\in(0,1]$, there exists a positive number $C_p$ such that
$$
\|Q_m\|_{L^p}\le C_pm^{1-1/p},\quad m\ge0.
$$

We are going to improve this estimate and show that one can select $C_p$ so that it remains bounded as $p$ approaches 1.

Denote by $\F$ and $\F^*$ the Fourier transform and the inverse Fourier transform defined by
$$
(\F f)(t)\df\int_\R f(x)e^{-2\pi\ri xt}\,dx,\quad(\F^*f)(t)\df\int_\R f(x)e^{2\pi\ri xt}\,dx.
$$

\begin{thm}
\label{ravnomerno}
Let $0<p\le1$ and let $q$ be an infinitely differentiable function on $\R$ with compact support.
Then the trigonometric polynomials $Q_m$ given by {\em\rf{Ups}} admit the estimate
$$
\|Q_m\|_{L^p(\T)}\le m^{1-1/p}\|\F q\|_{L^p(\R)}.
$$
\end{thm}

\Pf We define the operator $R$ from $L^p(\R)$ to $L^p(\T)$ by
$$
(Rf)\big(e^{2\pi\ri t}\big)=\sum_{k\in\Z}f(t+k).
$$
Clearly, for $f\in L^1(\R)$, we have
\bay
\label{Fourier}
\widehat{Rf}(n)=(\F f)(n),\quad n\in\Z.
\ey
It is easy to see that
$$
\|Rf\|^p_{L^p(\T)}\le\sum_{k\in\Z}\int_0^1|f(t+k)|^p\,dt=\|f\|_{L^p(\R)}^p.
$$

Put now 
$$
f=\F^*q\quad\mbox{and}\quad f_m(t)=mf(mt),\quad t\in\R.
$$
It follows from \rf{Fourier} that $(Rf_m)(t)=Q_m\big(e^{2\pi\ri t}\big)$, $t\in\R$. Thus,
$$
\|Q_m\|_{L^p(\T)}\le\left(\int_\R\big|mf(mt)|^p\,dt\right)^{1/p}
=m^{1-1/p}\|f\|_{L^p(\R)}=m^{1-1/p}\|\F q\|_{L^p(\R)}.\quad\bl
$$

\begin{cor}
Under the hypotheses of Theorem {\em{\ref{ravnomerno}}} for each $p_0\in(0,1)$,
\bay
\label{p0}
\sup_{p\in[p_0,1]}m^{1/p-1}\|Q_m\|_{L^p(\T)}<\be.
\ey
\end{cor}

We can improve now Lemma 4.2 of \cite{AP3} that shows to what extent the $L^p$-quasinorm of a trigonometric polynomial of a given degree can jump when we apply the Riesz projection. The improvement below 
gives us an estimate uniform in $p$ as $p$ approaches 1.

Recall that the Riesz projection 
$\pp_+$ and $\pp_-$ on the set of trigonometric polynomials are defined by
$$
(\pp_+\psi)(z)=\sum_{j\ge0} \widehat{\psi}(j)z^j\quad\mbox{and}
\quad(\pp_-\psi)(z)=\sum_{j<0} \widehat{\psi}(j)z^j.
$$

\begin{thm}
\label{rubim}
Let $q$ be an infinitely differentiable function on $\R$ such that 
$$
\supp q=[-1,1],\quad q(t)>0\quad\mbox{for}\quad t\in(-1,1)\quad\mbox{and}\quad q(0)=1.
$$ 
Consider the trigonometric polynomials  $Q_m$, $m\ge1$, of degree $m-1$ defined by {\em\rf{Ups}}. Then 
for each $p_0\in(0,1)$, there exists a positive number $d$ such that and a sequence of trigonometric polynomials $Q_m$, $m\ge1$, of degree $m-1$ such that
$$
\frac{\|\pp_+Q_m\|_{L^p}}{\|Q_m\|_{L^p}}\ge d m^{1/p-1}
$$
for every $p\ge p_0$.
\end{thm}

\medskip 

\Pf Clearly, 
$\deg Q_m=m-1$ and $(\pp_+Q_m)(0)=1$. It is also easy to see that
$$
\|\pp_+Q_m\|_{L^p}\ge|(\pp_+Q_m)(0)|=1.
$$
Denote the supremum of the left-hand side of \rf{p0} by $s$. Let $p\in[p_0,1)$.
Then by \rf{p0},
$$
\frac{\|\pp_+Q_m\|_{L^p}}{\|Q_m\|_{L^p}}\ge s^{-1} m^{1/p-1}.\quad\bl
$$

\medskip

%

\

\section{\bf Estimates of the $\bs{p}$-norms in $\bs{\bS_p}$ for certain Hankel matrices}
\setcounter{equation}{0}
\label{Gankeli}

\

In this section we recall certain facts on Hankel matrices (Hankel operators) that will be used in \S\;\ref{povedenie}. We also obtain sharp $\bS_p$ estimates for certain special Hankel matrices.

For a function $\f$ analytic in $\dd$, we denote by $\G_\f$ the {\it Hankel matrix}
$$
\G_\f=\{\widehat\f(j+k)\}_{j,k\ge0},
$$
where $\widehat\f(j)$ is the $j$th Taylor coefficient of $\f$. We refer the reader to the book \cite{Pe3} for basic definitions and properties of Hankel operators.

By Nehari's theorem, the matrix $\G_\f$ induces a bounded linear operator on the sequence space $\ell^2$ if and only if there exists a function $\psi$ in $L^\be(\T)$ such that their Fourier coefficients satisfy the equalities $\widehat\f(j)=\widehat\psi(j)$, $j\ge0$. Here $\widehat\psi(j)$ is the $j$th Fourier coefficient of $\psi$.

We will need the following criterion for the membership of $\G_\f$ in the Schatten--von Neumann class $\bS_p$: 

\medskip

{\it Let $0<p<\be$ and let $\f$ be a function analytic in $\dd$. Then
$$
\G_\f\in\bS_p\quad\Longleftrightarrow\quad\f\quad\mbox{belongs to the Besov class}\quad
B_p^{1/p}.
$$}

This was established in
\cite{Pe1} for $p\ge1$. For $p<1$, this was proved in \cite{Pe2}, see also \cite{Pek} and \cite{S} for alternative proofs. We also refer the reader to the book \cite{Pe3}, see Ch.6, \S\,3, in which a detailed presentation of this material is given.

Note that to prove the sufficiency of the condition $\f\in B_p^{1/p}$ 
for the membership of the Hankel operator $\G_\f$ in $\bS_p$
in the case $0<p\le1$, the following inequality was established in \cite{Pe1} and \cite{Pe2}, see also
the book \cite{Pe3}, Ch.6, \S\,3:
\bay
\label{dlya_polinomov}
\|\G_\f\|_{\bS_p}\le2^{1/p-1}m^{1/p}\|\f\|_{L^p}
\ey
for an arbitrary analytic polynomial $\f$ of degree at most $m-1$.

The following fact is an improvement of Theorem 3.1 of \cite{AP3}.

\begin{thm}
\label{spetsform}
Let $0<p_0<1$. Then there exists a positive number $d$ such that for $p\in[p_0,1)$ and
for a polynomial $\f$ is of the form
$$
\f(z)=\sum_{j=2^{n-1}+1}^{2^{n+1}-1}\widehat\f(j)z^j,
$$
the following inequalities hold:
\bay
\label{s_dvukh_storon}
d2^{(n+1)/p}\|\f\|_{L^p}\le\|\G_\f\|_{\bS_p}\le2^{(n+1)/p}\|\f\|_{L^p}.
\ey
\end{thm}

\Pf The upper estimate for $\|\G_\f\|_{\bS_p}$ in \rf{s_dvukh_storon} is an immediate consequence of 
\rf{dlya_polinomov}. 

To obtain the lower estimate for $\|\G_\f\|_{\bS_p}$, we observe that 
the reasoning given in the proof of Lemma 3.9 of Chapter 6 of the book \cite{Pe3} together with Theorem \ref{rubim}.
implies that
\bay
\label{f*Vk}
2^k\|\f*V_k\|_{L^p}^p\le\const\|\G_\f\|_{\bS_p}^p,\quad k\ge1.
\ey
It follows from the definition of $V_k$ given in \rf{Vn} that for $n\ge2$,
$$
\f=\f*V_{n-1}+\f*V_n+\f*V_{n+1}.
$$
Thus, we can conclude from \rf{f*Vk} that
$$
\|\f\|_{L^p}^p\le\sum_{j=-1}^1\|\f*V_{n+j}\|^p_{L^p}\le\const\cdot2^{n+1}\|\G_\f\|_{\bS_p}^p.
$$
This implies the lower estimate for $\|\G_\f\|_{\bS_p}$ in \rf{s_dvukh_storon}. $\bl$

\medskip

%


Let us proceed to Hankel Schur multipliers of $\bS_p$.
Let $0<p\le1$ and let $\f$ be an analytic polynomial of degree at most $m-1$. We need the following inequality:
\bay
\label{multnorm}
\|\G_\f\|_{\fM_p}\le(2m)^{1/p-1}\|\f\|_{L^p}.
\ey
It was established in \cite{Pe2}; the reader can also find the proof of this estimate in \cite{Pe3}, Ch.6, \S\,3.

\

\section{\bf The behaviour of the $\bs{p}$-norms in $\bs{H^p}$ of the analytic Dirichlet kernel}
\setcounter{equation}{0}
\label{analDir}

\

We define the analytic Dirichlet kernel $D_n^+$ by
$$
D_n^+(z)\df\sum_{j=0}^{n-1}z^j,\quad n\ge1.
$$
We will see in \S\:\ref{povedenie} that the $p$-norms $\|\cp_n\|_{\mB(\bS_p)}$ admit sharp lower and upper estimates in terms of the $p$-norms $\|D_n^+\|_{H^p}$. The following theorem gives us sharp upper and lower estimates for 
$\|D_n^+\|_{H^p}$.

\begin{thm} 
\label{Dn+Hp}
Let $\frac12\le p<1$. 
Then
\bay
\label{sdvukhstoron}
\frac{\sqrt 2(1-e^{-1})^2}{4\pi}\min\big\{(1-p)^{-1},\log n\big\}\le \|D^+_n\|_{H^p}\le
\min\big\{2(1-p)^{-1},\log 5n\big\}.
\ey
\end{thm}

\Pf Let us first proceed to the upper estimate for $\|D_n^+\|_{H^p}$.
We have
$$
|D^+_n(e^{it})|=\frac{|1-e^{int}|}{|1-e^{it}|}=\frac{|\sin nt/2|}{|\sin t/2|}\le\frac{\pi |\sin nt/2|}{|t|}
\quad\text{for}\quad t\in[-\pi,\pi].
$$
Let $p\in(0,1)$. Then
\bay
\label{1p}
\|D^+_n\|_{H^p}^p\le\pi^{p-1}\int_0^\pi\frac{|\sin nt/2|^p\, dt}{t^p}\le\pi^{p-1}\int_0^\pi\frac{dt}{t^p}=\frac1{1-p}
\ey
and
\bey
\|D^+_n\|_{H^p}\le\|D^+_n\|_{H^1}\le\int_0^\pi\frac{|\sin nt/2|\, dt}{t}=\int_0^{\pi n/2}\frac{|\sin x|\,dx}{x}\\
\le1+\int_1^{\pi n/2}\frac{dx}{x}=1+\log(\pi n/2)=\log(e\pi n/2)\le\log 5n.
\eey

Thus, we have proved that
$$
\|D^+_n\|_{H^p}^p\le\min((1-p)^{-1},(\log 5n)^p).
$$
Note that  $p\mapsto(1-p)^{1-1/p}$ is a a decreasing function on $(0,1)$,
$\lim_{p\to0}(1-p)^{1-1/p}=e$ and $\lim_{p\to1}(1-p)^{1-1/p}=1$.
Hence, \rf{1p} implies
$$
\|D^+_n\|_{H^p}\le(1-p)^{1-1/p}(1-p)^{-1}\le e(1-p)^{-1}.
$$
It is also clear that for $p\in\big[\frac12,1)$,
$$
\|D^+_n\|_{H^p}\le 2(1-p)^{-1}.
$$
Thus,
$$
\|D^+_n\|_{H^p}\le\min(2(1-p)^{-1},\log 5n).
$$


Let us obtain now the lower estimate for $\|D^+_n\|_{H^p}$. Clearly, $\|D^+_n\|_{H^p}\ge1$ for all $n\ge1$ and all $p>0$.
Besides, $\|D^+_1\|_{H^p}=1$ for all $p>0$.
Suppose now that $n\ge2$.

We have
$$
|D^+_n(e^{it})|=\frac{|1-e^{int}|}{|1-e^{it}|}=\frac{|\sin nt/2|}{|\sin t/2|}\ge\frac{2|\sin nt/2|}{|t|}
\quad\text{for}\quad t\in\R.
$$
Hence,
\bey
\|D^+_n\|_{H^p}^p\ge\frac{2^p}{\pi}\int_0^\pi\frac{|\sin nt/2|^p\, dt}{t^p}=\frac{2n^{p-1}}\pi\int_0^{\pi n/2}\frac{|\sin x|^p\,dx}{x^p}\\
=\frac{2n^{p-1}}\pi\sum_{k=1}^n\int_{\frac\pi 2(k-1)}^{\frac\pi 2k}\frac{|\sin x|^p\,dx}{x^p}.
\eey

To estimate the integral $\int\limits_{\frac\pi 2(k-1)}^{\frac\pi 2k}\dfrac{dx}{x^p}$ for given 
$k\in\Bbb N$, we put
$$
E_k\df\Big\{x\in\Big[\frac\pi 2(k-1),\frac\pi 2k\Big]:2\sin^2x\ge1\Big\}.
$$ 
Clearly, $E_k$ is an interval of length $\pi/4$ and
\bey
\int_{\frac\pi 2(k-1)}^{\frac\pi 2k}\frac{|\sin x|^p\,dx}{x^p}\ge2^{-p/2}\int_{E_k}\frac{dx}{x^p}
\ge2^{-p/2}\frac\pi 4\min_{x\in E_k}x^{-p}
\ge2^{-p/2}\frac\pi 4\Big(\frac\pi 2k\Big)^{-p}
=2^{p/2-2}\pi^{1-p}k^{-p}.
\eey
Thus,
\bey
\|D^+_n\|_{H^p}^p
\ge\frac{2^{p/2-1}n^{p-1}}{\pi^p}\sum_{k=1}^n k^{-p}
\ge \frac{2^{p/2-1}n^{p-1}}{\pi^p}\int_{1}^{n+1}\frac{dt}{t^p}\\
=\frac{2^{p/2-1}n^{p-1}}{\pi^p(1-p)}((n+1)^{1-p}-1)
\ge\frac{2^{p/2-1}(1-n^{p-1})}{\pi^p(1-p)}.
\eey

If $(1-p)\log n\ge1$, then
\bey
\|D^+_n\|_{H^p}^p\ge\frac{2^{p/2-1}(1-n^{p-1})}{\pi^p(1-p)}=
\frac{2^{p/2-1}(1-e^{(p-1)\log n})}{\pi^p(1-p)}
\ge\frac{2^{p/2-1}(1-e^{-1})}{\pi^p(1-p)}.
\eey

Thus, for $p\in[1/2,1)$ we have
\bey
\|D^+_n\|_{H^p}\ge\frac{2^{1/2-1/p}(1-e^{-1})^{1/p}}{\pi(1-p)^{1/p}}
\ge\frac{\sqrt 2(1-e^{-1})^2}{4\pi(1-p)^{1/p}}\ge\frac{\sqrt 2(1-e^{-1})^2}{4\pi(1-p)}
\eey
if $(1-p)\log n\ge1$.

Now let $(1-p)\log n\le1$. Then
\bey
\|D^+_n\|_{H^p}^p\ge\frac{2^{p/2-1}(1-n^{p-1})}{\pi^p(1-p)}=
\frac{2^{p/2-1}(1-e^{(p-1)\log n})}{\pi^p(1-p)\log n}\log n
\ge\frac{2^{p/2-1}(1-e^{-1})\log n}{\pi^p}.
\eey

Hence, for $p\in[1/2,1]$, we have
\bey
\|D^+_n\|_{H^p}\ge\frac{2^{1/2-1/p}(1-e^{-1})^{1/p}(\log n)^{1/p}}{\pi}\ge\frac{\sqrt 2(1-e^{-1})^2}{4\pi}\log n.
\quad\bl
\eey

\

\section{\bf The behaviour of $\bs{\|\cp_n\|_{\mB(\bS_p)}}$ versus the behavior of $\bs{\|D^+_n\|_{H^p}}$}
\setcounter{equation}{0}
\label{povedenie}

\

In this section we reduce the estimation of  $\|\cp_n\|_{\mB(\bS_p)}$ to the estimation of $\|D^+_n\|_{H^p}$.
The following theorem is the main result of the section.

\begin{thm}
\label{cherepDn+}
Let $\frac12\le p<1$. Then there exists a positive number $c$ such that
\bay
\label{sverkhu-snizu}
cn^{1/p-1}\|D^+_n\|_{H^p}\le\|\cp_n\|_{\mB(\bS_p)}=\|\chi_n\|_{\fM_p}=\|\D_n\|_{\fM_p}\le(2n)^{1/p-1}\|D^+_n\|_{H^p}.
\ey
\end{thm}

\Pf The upper estimate in \rf{sverkhu-snizu} for $\|\cp_n\|_{\mB(\bS_p)}$ for all $p\in(0,1)$ was mentioned in \cite{AP3} and follows immediately from \rf{multnorm}.

It suffices to establish the lower estimate for $n=2^{k}$. This is clear from \rf{sdvukhstoron}. 

Let us proceed to the lower estimate for $\|\cp_n\|_{\mB(\bS_p)}$ in \rf{sverkhu-snizu}. Let $q$ be a nonnegative infinitely differentiable even function on $\R$ such that 
$$
\supp q=[-1,1],\quad q(t)>0\quad\mbox{for}\quad t\in(-1,1)
$$
and
\bay
\label{q(t)+q(t-1)}
q(t)+q(t-1)=1\quad\mbox{for}\quad t\in[0,1].
\ey
Consider the trigonometric polynomials $Q_m$ defined by \rf{Ups}. 

It is more convenient for our purpose to consider instead of the Riesz projection $\pp_+$ the projection
$\pp^+$ defined on the class of trigonometric polynomials by
$$
(\pp^+\psi)(z)=\sum_{j\ge1} \widehat{\psi}(j)z^j.
$$

It follows immediately from \rf{q(t)+q(t-1)}
that
\bay
\label{Q+Q-Dn}
\bar z\pp^+Q_{m+1}+z^m\pp_-Q_{m+1}=D_m^+.
\ey
Put
$$
P_k(z)=z^{2^k}Q_{2^{k-1}}(z),\quad k\ge1,
$$
$$
P_k^-=z^{2^k}\pp_-Q_{2^{k-1}}\quad\mbox{and}\quad P_k^+=z^{2^k}\pp^+Q_{2^{k-1}}.
$$
It is easy to see that
$
P_k^-=P_k*D^+_{2^k}.
$
Thus,
$$
\G_{P_k^-}=\G_{P_k}\star\D_{2^k-1}.
$$

Let $m=2^{k-1}-1$. We have by \rf{Q+Q-Dn}
\bay
\label{Pk-+}
P_k^++z^{2^{k-1}}P_k^-=z^{2^k+1}D^+_{2^{k-1}-1}.
\ey
It follows from the fact that $q$ is an even function that
$$
\|P_k^+\|_{H^p}^p=\|P_k^-\|_{H^p}^p
$$
and so by \rf{Pk-+},
\bay
\label{D+Pk-}
\|D^+_{2^{k-1}-1}\|_{H^p}^p\le\|P_k^+\|_{H^p}^p+\|P_k^-\|_{H^p}^p=2\|P_k^-\|_{H^p}^p.
\ey

%



By Theorem \ref{spetsform},  there exists a positive number $d$ such that
inequalities \rf{s_dvukh_storon} hold for $\f=P_k$ and $\f=P_k^-$.

Thus,
$$
\|\D_{2^k}\|_{\fM_p}\ge\const\frac{\|\G_{P_k^-}\|_{\bS_p}}{\|\G_{P_k}\|_{\bS_p}}
\ge\const\frac{\|P_k^-\|_{H^p}}{\|P_k\|_{H^p}}.
$$
Applying \rf{D+Pk-}, we obtain
$$
\|\D_{2^k}\|_{\fM_p}\ge\const\frac{\|D^+_{2^{k-1}}\|_{H^p}}{\|P_k\|_{H^p}}
\ge\const2^{k(1/p-1)}\|D^+_{2^{k-1}}\|_{H^p}
$$
by Theorem \ref{ravnomerno}. It remains to observe that
$$
\|D^+_{2^{k-1}}\|_{H^p}\ge\const\|D^+_{2^{k}}\|_{H^p}
$$
which follows easily from Theorem \ref{Dn+Hp}.
This completes the proof. $\bl$


%
%
%
%


\

\section{\bf The main result}
\setcounter{equation}{0}
\label{Osnova}

\

We are ready now to state the main result of the paper.

\begin{thm}
\label{koda}
Let $\frac12\le p<1$. Then there exists a positive number $c$ such that
\begin{align*}
cn^{1/p-1}\min\big\{(1-p)^{-1},\log n\big\}&\le\|\cp_n\|_{\mB(\bS_p)}=
\|\chi_n\|_{\fM_p}=\|\D_n\|_{\fM_p}\\[.2cm]
&\le(2n)^{1/p-1}
\min\big\{2(1-p)^{-1},\log 5n\big\}.
\end{align*}
\end{thm}

Clearly, Theorem \ref{koda} is an immediate consequence of Theorems \ref{cherepDn+}
and \ref{Dn+Hp}.

Let now $n\ge3$ and $1-(\log n)^{-1}\le p<1$. Then there exist positive numbers $k$ and $K$ such that
$$
k\log n\le\|\cp_n\|_{\mB(\bS_p)}\le K\log n
$$
which implies the well known fact that 
$$
k\log n\le\|\cp_n\|_{\mB(\bS_1)}\le K\log n.
$$

\medskip

{\bf Remark.} We have observed in \cite{AP3} that if we consider the triangular projection $\cp_n^\K$ on the space of $n\times n$ block matrices whose entries are bounded linear operators on a Hilbert space $\K$, then
$$
\|\cp_n^\K\|_{\bS_p\to\bS_p}=\|\cp_n\|_{\bS_p\to\bS_p},
$$
and so in the statement of Theorem \ref{koda} we can replace $\|\cp_n\|_{\bS_p\to\bS_p}$ by 
$\|\cp_n^\K\|_{\bS_p\to\bS_p}$.

\

%
%
%
%
%
%

\

%
%
%
%
%
%

\
 
 \begin{footnotesize}
 
\noindent
\begin{tabular}{p{7cm}p{15cm}}
A.B. Aleksandrov & V.V. Peller \\
St.Petersburg Department & St.Petersburg State University\\
Steklov Institute of Mathematics  &Universitetskaya nab., 7/9 \\
Fontanka 27, 191023 St.Petersburg & 199034 St.Petersburg, Russia\\
Russia&\\
email: alex@pdmi.ras.ru&\\
\\
&St.Petersburg Department\\
&Steklov Institute of Mathematics\\
&Russian Academy of Sciences\\
&Fontanka 27, 191023 St.Petersburg\\
&Russia\\

& email: peller@math.msu.edu
\end{tabular}

\end{footnotesize}

\end{document}